\documentclass[12pt]{article}

\usepackage{amsmath}
\usepackage{amssymb}
\usepackage{mathdots}

\newtheorem{thm}{Theorem}[section]
\newtheorem{prop}[thm]{Proposition}
\newtheorem{cor}[thm]{Corollary}

\newtheorem{lem}[thm]{Lemma}

\newtheorem{ex}[thm]{Example}

\newcommand{\gp}{\mathfrak{p}}
\newcommand{\gq}{\mathfrak{q}}
\newcommand{\gm}{\mathfrak{m}}

\newcommand{\ga}{\mathfrak{a}}

\newcommand{\zz}{\mathbb{Z}}

\newcommand{\rees}[1]{{\mathrm R}(#1)}

\newcommand{\lc}[3]{{\mathrm H}^{#1}_{#2}(#3)} 

\newcommand{\length}[2]{\ell_{#1}(\,#2\,)}
\newcommand{\ass}[2]{{\rm Ass}_{#1}\,{#2}}
\newcommand{\assh}[2]{{\rm Assh}_{#1}\,#2}
\newcommand{\Min}[2]{{\rm Min}_{#1}\,{#2}}
\newcommand{\height}[2]{{\rm ht}_{#1}\,#2}
\newcommand{\dep}[2]{{\rm depth}_{#1}\,#2}
\newcommand{\grade}[1]{{\rm grade}\,#1}

\newcommand{\setmat}[3]{{\rm Mat}(#1, #2\,; #3)} 
\newcommand{\rdim}[1]{{\rm dim}\,{#1}} 
\newcommand{\sat}[1]{(#1)^{\rm sat}} 
\newcommand{\detid}[2]{{\rm I}_{#1}(#2)} 

\newcommand{\ra}{\longrightarrow}
\newcommand{\sra}{\rightarrow}

\newcommand{\shomom}[1]{\stackrel{#1}{\rightarrow}}

\begin{document}

\title{\Large 
Saturations of powers of certain determinantal ideals}

\author{
Kosuke Fukumuro, Taro Inagawa and Koji Nishida\footnote{
The last author is supported by KAKENHI (23540042)}
}
\date{
\small
Graduate School of Science, Chiba University, \\
1-33 Yayoi-Cho, Inage-Ku, Chiba-Shi, 263-8522, JAPAN
}
\maketitle

\begin{abstract}
Let $R$ be a Noetherian local ring and $m$ a positive integer.
Let $I$ be the ideal of $R$ generated by the maximal minors
of an $m \times (m + 1)$ matrix $M$ with entries in $R$\,.
Assuming that the grade of the ideal generated by 
the $k$-minors of $M$ is at least $m - k + 2$
for $1 \leq \forall k \leq m$\,,
we will study the associated primes of $I^n$ for $\forall n > 0$\,.
Moreover, we compute the saturation of $I^n$
for $1 \leq \forall n \leq m$ in the case where $R$ is
a Cohen-Macaulay ring
and the entries of $M$ are powers of elements
that form an sop for $R$\,.
\end{abstract}

\section{Introduction}

Let $(R, \gm)$ be a Noetherian local ring and
$I$ an ideal of $R$ such that $\rdim{R / I} > 0$\,.
Let $n$ be a positive integer.
We set
\[
\sat{I^n} \, = \, \bigcup_{i > 0}\,(I^n :_R \gm^i)
\]
and call it the saturation of $I^n$\,.
As $\sat{I^n} / I^n \cong \lc{0}{\gm}{R / I^n}$\,,
where $\lc{0}{\gm}{\,\cdot\,}$ denotes the $0$-th local cohomology functor,
we have $\sat{I^n} = I^n$ if and only if
$\dep{}{R / I^n} > 0$\,.
Moreover, if $J$ is an $\gm$-primary ideal such that
$\dep{}{R / (I^n :_R J)} > 0$\,,
we have $\sat{I^n} = I^n :_R J$\,.
On the other hand,
the $n$-th symbolic power of $I$ is defined by
\[
I^{(n)} = \bigcap_{\gp\,\in\,\Min{R}{R / I}}
(I^n R_\gp \cap R)\,.
\]
In order to compare $\sat{I^n}$ and $I^{(n)}$\,,
let us take a minimal primary decomposition of $I^n$\,;
\[
I^n = \bigcap_{\gp\,\in\,\ass{R}{R / I^n}} {\rm Q}(\gp)\,,
\]
where ${\rm Q}(\gp)$ denotes the $\gp$-primary component.
It is obvious that
\[
\sat{I^n} = \bigcap_{\gm\,\neq\,\gp\,\in\,\ass{R}{R / I^n}}{\rm Q}(\gp)
\hspace{3ex}\mbox{and}\hspace{3ex}
I^{(n)} = \bigcap_{\gp\,\in\,\Min{R}{R / I}}{\rm Q}(\gp)\,.
\]
Hence we have $\sat{I^n} \subseteq I^{(n)}$ and
the equality holds if and only if
$\ass{R}{R / I^n}$ is a subset of
$\{\,\gm\,\} \cup \Min{R}{R / I}$\,.
Therefore $\sat{I^n} = I^{(n)}$ if $\rdim{R / I} = 1$\,.
However, if $\rdim{R / I} \geq 2$\,,
usually it is not easy to decide whether
$\sat{I^n} = I^{(n)}$ or not.
Furthermore, describing a system of generators for
$\sat{I^n} / I^n$ precisely is often very hard.
In this paper,
assuming that $R$ is an $(m + 1)$-dimensional Cohen-Macaulay local ring
and $I$ is an ideal generated by the maximal minors of the following
$m \times (m + 1)$ matrix;
\[
M = \left(\begin{array}{@{\,}cccccc@{\,}}
x_1 & x_2 & x_3 & \cdots & x_m & x_{m + 1} \\
x_2 & x_3 & x_4 & \cdots & x_{m + 1} & x_1 \\
x_3 & x_4 & x_5 & \cdots & x_1 & x_2 \\
\hdotsfor{6} \\
x_m & x_{m + 1} & x_1 & \cdots & x_{m - 2} & x_{m - 1}
\end{array}\right)\,,
\]
where $x_1, x_2, \dots , x_{m + 1}$ is an sop for $R$\,,
we aim to prove
\begin{itemize}
\item
$\sat{I^n} \subsetneq I^{(n)}$ if $n > m \geq 3$\,,
\item
$\sat{I^n} = I^n$ if $1 \leq n < m$\,,
\item
$\sat{I^m} = I^m :_R (x_1, x_2, \dots , x_{m + 1})R$\,,
\item
$\sat{I^m} / I^m \cong R / (x_1, x_2, \dots , x_{m + 1})R$ if $m \geq 2$\,.
\end{itemize}
Moreover, we describe a generator of 
$\sat{I^m} / I^m$ using the determinant of a certain matrix induced from $M$\,.
The proofs of the assertions stated above are given in
Section 3 and 4 taking more general matrices as $M$\,.

Throughout this paper $R$ is a commutative ring,
and we often assume that $R$ is a Noetherian local ring
with the maximal ideal $\gm$\,.
For positive integers $m, n$ and an ideal $\ga$ of $R$\,, we denote by
$\setmat{m}{n}{\ga}$ the set of $m \times n$ matrices with entries in $\ga$\,.
For any $M \in \setmat{m}{n}{R}$ and any $k \in \zz$
we denote by $\detid{k}{M}$ the ideal generated by the
$k$-minors of $M$\,.
In particular, $\detid{k}{M}$ is defined to be $R$
(resp. $(0)$) for $k \leq 0$ (resp. $k > \min \{ m, n \}$).
If $M, N \in \setmat{m}{n}{R}$ and the $(i, j)$ entries of $M$ and $N$
are congruent modulo a fixed ideal $\ga$ for $\forall (i, j)$\,,
we write $M \equiv N \,\mbox{mod}\, \ga$\,.

\section{Preliminaries}

In this section, we assume that
$R$ is just a commutative ring.
Let $m, n$ be positive integers with $m \leq n$
and $M = (\,x_{ij}\,) \in \setmat{m}{n}{R}$\,.
Let us recall the following rather well-known fact.

\begin{lem}\label{2.1}
Suppose $\detid{m}{M} \subseteq \gp \in {\rm Spec}\,R$ and put
$\ell = \max\{\,0 \leq k \in \zz \mid
\detid{k}{M} \not\subseteq \gp\,\}$\,.
Then $\ell < m$ and there exists
$N \in \setmat{m - \ell}{n - \ell}{\gp R_\gp}$ such that
$\detid{k}{M}_\gp = \detid{k - \ell}{N}$
for any $k \in \zz$\,.
\end{lem}

\noindent
{\it Proof}. \,
We prove by induction on $\ell$\,.
The assertion is obvious if $\ell = 0$\,.
So, let us consider the case where $\ell > 0$\,.
Then $\detid{1}{M} \not\subseteq \gp$\,,
and so some entry of $M$ is a unit in $R_\gp$\,.
Hence, applying elementary operations to $M$
in $\setmat{m}{n}{R_\gp}$\,,
we get a matrix of the form
\[
\left(\begin{array}{@{\,}c|ccc@{\,}}
1 & 0 & \cdots & 0 \\ \hline
0 &   &   &   \\
\vdots & \multicolumn{3}{c}{\raisebox{0pt}[0pt][0pt]{\large $M'$}} \\
0 &   &   & 
\end{array}\right) \,,
\]
where $M' \in \setmat{m-1}{n-1}{R_\gp}$\,.
It is easy to see that
$\detid{k}{M}_\gp = \detid{k-1}{M'}$
for any $k \in \zz$\,. Hence
$\ell - 1 =
\max\{\,0 \leq k \leq m - 1 \mid
\detid{k}{M'} \not\subseteq \gp R_\gp \,\}$\,.
By the hypothesis of induction,
there exists
\[
N \in \setmat{\,(m-1)-(\ell - 1)\,}{\,(n-1)-(\ell - 1)\,}{\gp R_\gp} =
\setmat{m-\ell}{n-\ell}{\gp R_\gp}
\]
such that
$\detid{t}{M'} = \detid{t-(\ell -1)}{N}$ for any $t \in \zz$\,.
Then we have
$\detid{k}{M}_\gp = \detid{k-\ell}{N}$ for any $k \in \zz$\,.

\vspace{1em}
In the rest of this section,
we assume $n = m + 1$\,.
For $1 \leq \forall j \leq m + 1$\,,
$M_j$ denotes the $m \times m$ submatrix of $M$
determined by removing the $j$-th column.
We set $a_j = (-1)^{j-1}\cdot\det M_j$ and
$I = (a_1, a_2, \dots , a_{m+1})R = \detid{m}{M}$\,.
Let us take an indeterminate $t$ over $R$ and consider
the Rees algebra of $I$\,;
\[
\rees{I}:= R[a_1t, a_2t , \dots , a_{m+1}t] \subseteq R[\,t\,]\,,
\]
which is a graded ring such that $\deg a_jt = 1$
for $1 \leq \forall j \leq m+1$\,.
On the other hand,
let $S = R[T_1, T_2, \dots , T_{m+1}]$ be a polynomial ring
over $R$ with $m+1$ variables.
We regard $S$ as a graded ring by setting $\deg T_j = 1$
for $1 \leq \forall j \leq m+1$\,.
Let $\pi : S \ra \rees{I}$ be the homomorphism of $R$-algebras
such that $\pi(T_j) = a_jt$ for $1 \leq \forall j \leq m+1$\,.
Then $\pi$ is a surjective graded homomorphism.
Now we set
\[
f_i = \sum_{j = 1}^{m+1}\,x_{ij}T_j \in S_1
\]
for $1 \leq \forall i \leq m$\,.
It is easy to see $(f_1, f_2 , \dots , f_m)S \subseteq {\rm Ker}\,\pi$\,.
For our purpose, the following result due to Avramov \cite{a} is very important
(Another elementary proof is given in \cite{f}).

\begin{thm}\label{2.2}
Suppose that $R$ is a Noetherian ring.
If $\grade{\detid{k}{M}} \geq m - k + 2$ for $1 \leq \forall k \leq m$\,,
then ${\rm Ker}\,\pi = (f_1, f_2, \dots , f_m)S$ and
$\grade{(f_1, f_2, \dots , f_m)S} = m$\,.
\end{thm}

As the last preliminary result,
we describe a technique using determinants of matrices.
Suppose that $y_1, y_2, \dots , y_{m+1}$ are elements of $R$ such that
\[
M
\left(\begin{array}{c}
y_1 \\ y_2 \\ \vdots \\ y_{m+1}
\end{array}\right) =
\left(\begin{array}{c}
0 \\ 0 \\ \vdots \\ 0
\end{array}\right)\,.
\]
We put $y = y_1 + y_2 + \cdots + y_{m+1}$\,.

\begin{lem}\label{2.3}
If $y, y_k$ form a regular sequence for some $k$
with $1 \leq k \leq m + 1$\,,
then there exists $\delta \in R$ such that
$y_j\cdot\delta = a_j$ for $1 \leq \forall j \leq m + 1$\,.
\end{lem}

\noindent
{\it Proof.}\,
We put $a = a_1 + a_2 + \cdots + a_{m+1}$\,.
Then the following assertion holds:

\vspace{1em}
\noindent
{\bf Claim}\,
$y\cdot a_j = y_j\cdot a$ for $1 \leq \forall j \leq m + 1$\,.

\vspace{1em}
\noindent
In order to prove the claim above,
let us consider the following
$(m + 1) \times (m + 1)$ matrix:
\[
N = \left(\begin{array}{cccc}
1 & 1 & \cdots & 1 \\ \hline
   &   &    &    \\
\multicolumn{4}{c}{\raisebox{0pt}[0pt][0pt]{\large $M$}}  \\
   &   &    &
\end{array}\right)\,.
\]
Expanding $\det N$ along the first row,
we get $\det N = a$\,.
We fix $j$ with $1 \leq j \leq m + 1$\,.
Multiplying the $j$-th column of $N$ by $y_j$\,,
we get
\[
N' =
\bordermatrix{
  &  &  & \stackrel{j}{\smallsmile} &  &   \cr
  & 1 & \cdots & y_j & \cdots & 1  \cr
  & x_{11} & \cdots & x_{1j}y_j & \cdots & x_{1, m+1} \cr
  & \vdots &   & \vdots &    & \vdots  \cr
  & x_{m1} & \cdots & x_{mj}y_j & \cdots & x_{m, m+1}
}\,.
\]
Then $\det N' = y_j\cdot\det N = y_j \cdot a$\,.
Next, for $1 \leq \forall \ell \leq m + 1$ with
$\ell \neq j$\,,
we add the $\ell$-th column of $N'$ multiplied by $y_\ell$ to
the $j$-th column, and get
\[
N'' =
\bordermatrix{
  &  &  & \stackrel{j}{\smallsmile} &  &   \cr
  & 1 & \cdots & y & \cdots & 1  \cr
  & x_{11} & \cdots & 0 & \cdots & x_{1, m+1} \cr
  & \vdots &   & \vdots &    & \vdots  \cr
  & x_{m1} & \cdots & 0 & \cdots & x_{m, m+1}
}\,,
\]
since our assumption means
\[
x_{i1}y_1 + \cdots + x_{ij}y_j + \cdots + x_{i, m+1}y_{m+1} = 0
\]
for $1 \leq \forall i \leq m$\,.
Then $\det N'' = \det N' = y_j\cdot a$\,.
Finally, replacing the first $j$ columns of $N''$\,, we get
\[
N''' =
\left(
\begin{array}{@{\,}c|ccc@{\,}}
y & 1 & \cdots & 1 \\ \hline
0 &   &   &   \\
\vdots & \multicolumn{3}{c}{\raisebox{0pt}[0pt][0pt]{\large $M_j$}} \\
0 &   &   & 
\end{array}
\right) \,.
\]
Then $y\cdot a_j = y\cdot (-1)^{j-1}\cdot\det M_j =
(-1)^{j-1}\cdot\det N''' = \det N'' = y_j\cdot a$\,.
Thus we get the equalities of the claim.

Now we take $k$ with $1 \leq k \leq m + 1$ so that
$y, y_k$ form a regular sequence.
Because $y\cdot a_k = y_k\cdot a$\,,
there exists $\delta \in R$ such that $a = y\delta$\,.
Then $y\cdot a_j = y_j\cdot y\delta$ for $1 \leq \forall j \leq m + 1$\,.
As $y$ is an $R$-NZD,
we get $a_j = y_j\cdot \delta$ for $1 \leq \forall j \leq m + 1$\,,
and the proof is complete.

\begin{lem}\label{2.4}
If $R$ is a Cohen-Macaulay local ring and
$y_1, y_2, \dots , y_{m+1}$ is an ssop for $R$\,,
then $y, y_k$ form a regular sequence for $1 \leq \forall k \leq m + 1$\,.
\end{lem}

\noindent
{\it Proof.}\,
It is enough to show for $k = 1$\,.
Because $(y_1, y_2, \dots , y_{m+1})R = (y, y_1, \dots , y_m)R$\,,
it follows that $y, y_1, \dots , y_m$ is an ssop for $R$\,, too.
Hence $y, y_1$ is $R$-regular.

\begin{lem}\label{2.5}
Suppose that $\ga$ is an ideal of $R$ and
$x_{ij} \in \ga$ for $\forall i, \forall j$\,.
We put $Q = (y_1, y_2, \dots , y_{m+1})R$\,.
Then $\delta$ of \ref{2.3} is an element of $\ga^m :_R Q$\,.
\end{lem}

\noindent
{\it Proof.}\,
We get this assertion since
$a_j \in \ga^m$ for $1 \leq \forall j \leq m + 1$\,.

\section{Associated primes of $R / I^n$}

Let $R$ be a Noetherian ring and
$M = ( x_{ij} ) \in \setmat{m}{m+1}{R}$\,,
where $1 \leq m \in \zz$\,.
Let $I = \detid{m}{M}$\,.
Throughout this section,
we assume that $I$ is a proper ideal and
$\grade{\detid{k}{M}} \geq m - k + 2$ for $1 \leq \forall k \leq m$\,.
Let us keep the notations of Section 2.

Let $K_\bullet$ be the Koszul complex of $f_1, f_2, \dots , f_m$\,,
which is a complex of graded free $S$-modules.
We denote its boundary map by $\partial_\bullet$\,.
Let $e_1, e_2, \dots , e_m$ be an $S$-free basis of $K_1$
consisting of homogeneous elements of degree $1$ such that
$\partial_1 (e_i) = f_i$ for $1 \leq \forall i \leq m$\,.
Then, for $1 \leq \forall r \leq m$\,,
\[
\{\,e_{i_1} \wedge e_{i_2} \wedge \cdots \wedge e_{i_r} \mid
1 \leq i_1 < i_2 < \cdots < i_r \leq m\,\}
\]
is an $S$-free basis of $K_r$ consisting of homogeneous elements
of degree $r$, and we have
\[
\partial_r ( e_{i_1} \wedge e_{i_2} \wedge \cdots \wedge e_{i_r} ) =
\sum_{p = 1}^r  (-1)^{p-1}\cdot f_{i_p}\cdot
e_{i_1} \wedge \cdots \wedge \widehat{e_{i_p}} \wedge 
\cdots \wedge e_{i_r}\,.
\]
Let $1 \leq n \in \zz$\,.
Taking the homogeneous part of degree $n$ of $K_\bullet$\,,
we get a complex
\[
[K_\bullet]_n \,:\,
0 \ra [K_m]_n \stackrel{\partial_m}{\ra} [K_{m-1}]_n \ra \cdots \ra
[K_1]_n \stackrel{\partial_1}{\ra} [K_0]_n \ra 0
\]
of finitely generated free $R$-modules.
It is obvious that $[K_r]_n = 0$ if $n < r$\,.
On the other hand, if $n \geq r$\,, then
\[
\left\{
T_1^{\,\alpha_1}T_2^{\,\alpha_2}\cdots T_{m+1}^{\,\alpha_{m+1}}\cdot
e_{i_1} \wedge e_{i_2} \wedge \cdots \wedge e_{i_r}
\hspace{1ex}\left| \hspace{1ex}
\begin{array}{l}
0 \leq \alpha_1, \alpha_2, \dots , \alpha_{m+1} \in \zz\,, \\
\alpha_1 + \alpha_2 + \cdots + \alpha_{m+1} = n - r\,, \\
1 \leq i_1 < i_2 < \cdots < i_r \leq m
\end{array}
\right.
\right\}
\]
is an $R$-free basis of $[K_r]_n$\,.

\begin{prop}\label{3.1}
If $(R , \gm)$ is a local ring
and $M \in \setmat{m}{m+1}{\gm}$\,,
we have
\[
{\rm proj.\,dim}_R\,R / I^n = 
\left\{
\begin{array}{ll}
n + 1 & \mbox{if $n < m$}\,,  \vspace{0.5em}\\
m + 1 & \mbox{if $n \geq m$}\,.
\end{array}\right.
\]
\end{prop}

\noindent
{\it Proof.}\,
By \ref{2.2} and \cite[1.6.17]{bh}, we see that
\[
0 \ra K_m \stackrel{\partial_m}{\ra} K_{m-1} \ra \cdots \ra
K_1 \stackrel{\partial_1}{\ra} K_0 \stackrel{\pi}{\ra}
\rees{I} \ra 0
\]
is a graded $S$-free resolution of $\rees{I}$\,.
Hence, for $0 \leq \forall n \in \zz$\,,
\[
0 \ra [K_m]_n \stackrel{\partial_m}{\ra} [K_{m-1}]_n \ra \cdots \ra
[K_1]_n \stackrel{\partial_1}{\ra} [K_0]_n \stackrel{\pi}{\ra}
I^nt^n \ra 0
\]
is an $R$-free resolution of the $R$-module $I^nt^n$\,.
Let us notice $I^nt^n \cong I^n$ as $R$-modules.
Suppose $1 \leq r \leq m$ and $n \geq r$\,.
Then, for any non-negative integers $\alpha_1, \alpha_2, \dots , \alpha_{m+1}$
with $\alpha_1 + \alpha_2 + \cdots + \alpha_{m+1} = n - r$ and
positive integers $i_1, i_2, \dots , i_r$ with
$1 \leq  i_1 < i_2 < \cdots < i_r \leq m$\,, we have
\begin{eqnarray*}
 &    & \partial_r( T_1^{\,\alpha_1}T_2^{\,\alpha_2}\cdots T_{m+1}^{\,\alpha_{m+1}}\cdot
                e_{i_1} \wedge e_{i_2} \wedge \cdots \wedge e_{i_r} )   \\
 & = & T_1^{\,\alpha_1}T_2^{\,\alpha_2}\cdots T_{m+1}^{\,\alpha_{m+1}}\cdot
              \sum_{p = 1}^r (-1)^{p-1}\cdot (\,\sum_{j = 1}^{m+1} x_{i_p, j}T_j\,)\cdot
              e_{i_1} \wedge \cdots \wedge \widehat{e_{i_p}} \wedge 
              \cdots \wedge e_{i_r}    \\
 & = & \sum_{p = 1}^r \sum_{j = 1}^{m+1} (-1)^{p-1}x_{i_p, j}\cdot
              T_1^{\,\alpha_1}\cdots T_j^{\,1 + \alpha_j}\cdots T_{m+1}^{\,\alpha_{m+1}}\cdot
              e_{i_1} \wedge \cdots \wedge \widehat{e_{i_p}} \wedge 
              \cdots \wedge e_{i_r} \\
 & \in & \gm\cdot [K_{r-1}]_n\,.
\end{eqnarray*}
Hence $[K_\bullet]_n$ gives a minimal $R$-free resolution of $I^n$\,.
If $n < m$\,, we have $[K_n]_n \neq 0$ and
$[K_r]_n = 0$ for any $r > n$\,,
and so ${\rm proj.\,dim}_R\,I^n = n$\,.
On the other hand,
if $n \geq m$\,,
we have $[K_m]_n \neq 0$ and $[K_r]_n = 0$ for any $r > m$\,,
and so ${\rm proj.\,dim}_R\,I^n = m$\,.
Thus we get the required equality as
${\rm proj.\,dom}_R\,R / I^n = {\rm proj.\,dim}_R\,I^n + 1$\,. 

\vspace{1em}
By Auslander-Buchsbaum formula (cf. \cite[1.3.3]{bh}),
we get the following.

\begin{cor}\label{3.2}
If $(R, \gm)$ is local and $M \in \setmat{m}{m+1}{\gm}$, we have
\[
\dep{}{R / I^n} = \left\{
\begin{array}{ll}
\dep{}{R} - n - 1 & \mbox{if $n < m$}\,, \vspace{0.5em} \\
\dep{}{R} - m - 1 & \mbox{if $n \geq m$}\,.
\end{array}
\right.
\]
\end{cor}

\noindent
Here we remark that
$\dep{}{R} \geq \grade{\detid{1}{M}} \geq m + 1$
by our assumption of this section.
As a consequence of \ref{3.2},
we see that the next assertion holds.

\begin{cor}\label{3.3}
Suppose that $(R, \gm)$ is a local ring and
$M \in \setmat{m}{m+1}{\gm}$\,.
Then we have 
$\gm \in \ass{R}{R / I^n}$ if and only if
$n \geq m$ and $\dep{}{R} = m + 1$\,.
\end{cor}

The next result is a generalization of \ref{3.3}\,.

\begin{thm}\label{3.4}
Let $I \subseteq \gp \in {\rm Spec}\,R$ and $1 \leq n \in \zz$\,.
We put
$\ell = \max \{\,0 \leq k < m \mid \detid{k}{M} \not\subseteq \gp\,\}$\,.
Then the following conditions are equivalent.
\begin{itemize}
\item[{\rm (1)}]
$\gp \in \ass{R}{R / I^n}$\,.
\item[{\rm (2)}]
$n \geq m - \ell$ and $\dep{}{R_\gp} = m - \ell + 1$\,.
\end{itemize}
When this is the case, $\grade{\detid{\ell+1}{M}} = m - \ell + 1$\,.
\end{thm}

\noindent
{\it Proof.}\,
By \ref{2.1}, there exists
$N \in \setmat{m - \ell}{m - \ell + 1}{\gp R_\gp}$
such that $\detid{k}{N} = \detid{k+\ell}{M}_\gp$
for any $k \in \zz$\,.
Hence, for $1 \leq \forall k \leq m - \ell$\,,
we have
\[
\grade{\detid{k}{N}} = \grade{\detid{k+\ell}{M}_\gp} \geq
m - (k + \ell) + 2 = (m - \ell) - k + 2\,.
\]
Therefore, by \ref{3.3}, we see that
$\gp R_\gp \in \ass{R_\gp}{R_\gp / {\detid{m - \ell}{N}}^n}$
if and only if
$n \geq m - \ell$ and $\dep{}{R_\gp} = m - \ell + 1$\,.
Let us notice $\detid{m-\ell}{N} = I_\gp$\,.
Because $\gp \in \ass{R}{R / I^n}$ if and only if
$\gp R_\gp \in \ass{R_\gp}{R_\gp / {I_\gp}^n}$\,,
we see (1) $\Leftrightarrow$ (2)\,.
Furthermore, as $\detid{\ell + 1}{M} \subseteq \gp$\,,
we have $\grade{\detid{\ell + 1}{M}} \leq \dep{}{R_\gp}$\,,
and so we get $\grade{\detid{\ell + 1}{M}} = m - \ell + 1$
if the condition (2) is satisfied.

\vspace{1em}
For $1 \leq \forall n \in \zz$\,,
let $\Lambda^n_M$ be the set of integers $i$ such that
$\max\{  1, m - n + 1 \} \leq i \leq m$ and
$\grade{\detid{i}{M}} = m - i + 2$\,.
Then the following assertion holds.

\begin{thm}\label{3.5}
Let $R$ be a Cohen-Macaulay ring.
Then, for $1 \leq \forall n \in \zz$\,, we have
\[
\ass{R}{R / I^n} = \bigcup_{i \in \Lambda^n_M} 
\assh{R}{R / \detid{i}{M}}\,.
\]
\end{thm}

\noindent
{\it Proof.}\,
Let us take any $\gp \in \ass{R}{R / I^n}$ and put
$\ell = \max\{ 0 \leq k < m \mid \detid{k}{M} \not\subseteq \gp \}$\,.
Then $\detid{\ell + 1}{M} \subseteq \gp$\,.
Moreover, by \ref{3.4} we have $n \geq m - \ell$\,,
$\dep{}{R_\gp} = m - \ell + 1$ and
$\grade{\detid{\ell + 1}{M}} = m - \ell + 1$\,.
Hence $\ell + 1 \in \Lambda^n_M$\,.
Let us notice that
$\height{}{\gp} = \dep{}{R_\gp}$ and
$\height{}{\detid{\ell + 1}{M}} = \grade{\detid{\ell + 1}{M}}$ as
$R$ is Cohen-Macaulay.
Therefore $\height{}{\gp} = \height{}{\detid{\ell + 1}{M}}$\,,
which means
$\gp \in \assh{R}{R / \detid{\ell + 1}{M}}$\,.

Conversely, let us take any $i \in \Lambda^n_M$ and
$\gq \in \assh{R}{R / \detid{i}{M}}$\,.
Then $\height{}{\gq} = \height{}{\detid{i}{M}} =
\grade{\detid{i}{M}} = m - i + 2$\,.
As our assumption implies
$\height{}{\detid{i - 1}{M}} \geq m - i + 3$\,,
it follows that
$i - 1 = \max\{ 0 \leq k < m \mid \detid{k}{M} \not\subseteq \gq \}$\,.
Let us notice $n \geq m - (i - 1)$ as $m - n + 1 \leq i$\,,
which is one of the conditions for $i \in \Lambda^n_M$\,.
Moreover, we have $\dep{}{R_\gq} = \height{}{\gq} = m - (i - 1) + 1$\,.
Thus we get $\gq \in \ass{R}{R / I^n}$ by \ref{3.4}, and the proof is complete.

\begin{ex}\label{3.6}
Let $1 \leq m \in \zz$ and
let $R$ be an $(m + 1)$-dimensional Cohen-Macaulay local ring
with the maximal ideal $\gm$\,.
We take an sop
$x_1, x_2, \dots , x_{m + 1}$ for $R$ and a family
$\{ \alpha_{ij} \}_{1 \leq i \leq m, 1 \leq j \leq m+1}$
of positive integers.
For $1 \leq \forall i \leq m$ and $1 \leq \forall j \leq m + 1$\,, we set
\[
x_{ij} = \left\{
\begin{array}{ll}
x_{i + j - 1}^{\,\alpha_{ij}} & \mbox{if $i + j \leq m + 2$}\,,  \vspace{1em}\\
x_{i + j - m - 2}^{\,\alpha_{ij}} & \mbox{if $i + j > m + 2$}\,.
\end{array}
\right.
\]
Let us consider the following matrix;
\[
M = ( x_{ij} ) = 
\left(
\begin{array}{@{\,}cccccc@{\,}}
x_1^{\,\alpha_{11}} & x_2^{\,\alpha_{12}} & x_3^{\,\alpha_{13}}
      & \cdots & x_m^{\,\alpha_{1m}} & x_{m + 1}^{\,\alpha_{1, m+1}} \\
x_2^{\,\alpha_{21}} & x_3^{\,\alpha_{22}} & x_4^{\,\alpha_{23}}
      & \cdots & x_{m + 1}^{\,\alpha_{2m}} & x_1^{\,\alpha_{ 2, m+1}}  \\
x_3^{\,\alpha_{31}} & x_4^{\,\alpha_{32}} & x_5^{\,\alpha_{33}}
      & \cdots & x_1^{\,\alpha_{3m}} & x_2^{\,\alpha_{3, m+1}} \\
\hdotsfor{6} \\
x_m^{\,\alpha_{m1}} & x_{m + 1}^{\,\alpha_{m2}} & x_1^{\,\alpha_{m3}}
      & \cdots & x_{m - 2}^{\,\alpha_{mm}} & x_{m - 1}^{\,\alpha_{m, m+1}}
\end{array}\right)\,.
\]
If $\alpha_{ij} = 1$ for $\forall i$ and $\forall j$\,,
then $M$ is the matrix stated in Introduction.
We put $I = \detid{m}{M}$\,.
Then the following assertions hold.
\begin{itemize}
\item[{\rm (1)}]
$\height{}{\detid{k}{M}} \geq m - k + 2$ for $1 \leq \forall k \leq m$\,.
\item[{\rm (2)}]
${\displaystyle
{\rm proj.\,dim}_R\,R / I^n = \left\{
\begin{array}{ll}
n + 1 & \mbox{if $n < m$}\,, \vspace{1em}\\
m + 1 & \mbox{if $n \geq m$}\,.
\end{array}
\right.
}$
\item[{\rm (3)}]
${\displaystyle
\dep{}{R / I^n} = \left\{
\begin{array}{ll}
m - n & \mbox{if $n < m$}\,, \vspace{1em}\\
0 & \mbox{if $n \geq m$}\,.
\end{array}
\right.
}$
\end{itemize}
Furthermore, if $\alpha_{ij} = 1$ for $\forall i$ and $\forall j$\,,
the following assertions hold.
\begin{itemize}
\item[{\rm (4)}]
If $m \geq 2$\,,
then $\height{}{\detid{2}{M}} = m$ and
$\assh{R}{R / \detid{2}{M}} \subseteq \ass{R}{R / I^n}$
for any $n \geq m - 1$\,.
\item[{\rm (5)}]
If $m$ is an odd integer with $m \geq 3$\,, then
$\height{}{\detid{3}{M}} = m - 1$ and
$\assh{R}{R / \detid{3}{M}} \subseteq \ass{R}{R / I^n}$
for any $n \geq m - 2$\,.
\item[{\rm (6)}]
If $m \geq 3$\,, then
$\sat{I^n} \subsetneq I^{(n)}$ for any $n \geq m - 1$\,.
\end{itemize}
\end{ex}

\noindent
{\it Proof.}\,
(1)\,
We aim to prove the following.

\vspace{1em}
\noindent
{\bf Claim}\,
$J_{k-1} + \detid{k}{M}$ is $\gm$-primary for
$1 \leq \forall k \leq m + 1$\,,
where $J_{k-1} = (x_1, x_2, \dots , x_{k-1})R$\,.

\vspace{1em}
\noindent
If this is true, we have
$\dim R / \detid{k}{M} \leq k - 1$\,, and so
$\height{}{\detid{k}{M}} \geq \dim R - (k - 1) = m - k + 2$\,,
which is the required inequality.

In order to prove Claim, we take any $\gp \in {\rm Spec}\,R$
containing $J_{k-1} + \detid{k}{M}$\,.
It is enough to show $J_{m+1} 
= (x_1, x_2, \dots , x_{m+1})R \subseteq \gp$\,.
For that purpose, we prove $J_{\ell} \subseteq \gp$
for $k - 1 \leq \forall \ell \leq m + 1$ by induction on $\ell$\,.
As we obviously have $J_{k-1} \subseteq \gp$\,,
let us assume $k \leq \ell \leq m + 1$ and $J_{\ell - 1} \subseteq \gp$\,.
Because the $k$-minor of $M$ with respect to the first $k$ rows
and the columns $\ell - k + 1, \dots , \ell - 1, \ell$ is congruent with
\[
\det
\left(
\begin{array}{cccc}
   &    &    & x_\ell^{\,\alpha_{1, \ell}} \\
\multicolumn{2}{c}{\raisebox{0pt}[0pt][0pt]{\LARGE $0$}} 
            & x_\ell^{\,\alpha_{2, \ell -1}} &    \\
   & \iddots & \multicolumn{2}{c}{\raisebox{0pt}[0pt][0pt]{\LARGE $\ast$}} \\
x_\ell^{\,\alpha_{k, \ell - k+1}} &   &   & 
\end{array}\right)
\]
mod $J_{\ell - 1}$\,,
it follows that $J_{\ell - 1} + \detid{k}{M}$ includes some power of $x_\ell$\,.
Hence $x_\ell \in \gp$\,,
and so we get $J_\ell \subseteq \gp$\,.

(2) and (3) follow from \ref{3.1} and \ref{3.2}, respectively.

In the rest of this proof, we assume
$\alpha_{ij} = 1$ for any $\forall i$ and $\forall j$\,.

(4)\,
Let $\gq = (x_1 - x_2, x_2 - x_3, \dots , x_m - x_{m + 1})R$\,.
Then $x_1 \equiv x_i$ mod $\gq$ for $1 \leq \forall i \leq m + 1$\,.
Hence, any $2$-minor of $M$ is congruent with
\[
\det
\left(
\begin{array}{cc}
x_1 & x_1 \\
x_1 & x_1
\end{array}
\right)
= 0
\]
mod $\gq$\,.
This means $\detid{2}{M} \subseteq \gq$\,,
and so $\height{}{\detid{2}{M}} \leq \mu_{R}(\gq) = m$\,.
On the other hand, 
$\height{}{\detid{2}{M}} \geq m$ by (1).
Thus we get $\height{}{\detid{2}{M}} = m$\,.
Then, for any $n \geq m - 1$\,,
we have $2 \in \Lambda^n_M$\,,
and so $\assh{R}{R / \detid{2}{M}} \subseteq \ass{R}{R / I^n}$
by \ref{3.5}.

(5)\,
Let $\gp$ be the ideal of $R$ generated by
$\{ x_i - x_{i + 2} \}$\,, where
$i$ runs all odd integers with $2 \leq i \leq m - 2$\,.
Similarly, we set $\gq$ to be the ideal of $R$ generated by
$\{ x_j - x_{j + 2} \}$\,,
where $j$ runs all even integers with $2 \leq j \leq m - 1$\,.
Let $M'$ be the submatrix of $M$ with the rows $i_1, i_2, i_3$
and the columns $j_1, j_2, j_3$\,, where
$1 \leq i_1 < i_2 < i_3 \leq m$ and
$1 \leq j_1 < j_2 < j_3 \leq m + 1$\,.
We can choose $p, q$ with $1 \leq p < q \leq 3$ so that
$i_p \equiv i_q$ mod $2$\,.
Then, for $1 \leq \forall r \leq 3$\,,
we have $i_p + j_r \equiv i_q + j_r$ mod $2$\,,
and so, if $i_p + j_r$ is odd (resp. even),
it follows that $x_{i_p, j_r} \equiv x_{i_q, j_r}$ mod $\gp$
(resp. $\gq$).
Hence, we see that the $p$-th row of $M'$ is congruent with
the $q$-th row of $M'$ mod $\gp + \gq$\,,
which means $\det M' \equiv 0$ mod $\gp + \gq$\,.
As a consequence, we get
$\detid{3}{M} \subseteq \gp + \gq$\,.
Therefore $\height{}{\detid{3}{M}} \leq 
\mu_R(\gp) + \mu_R(\gq) = (m - 1) / 2 + (m - 1) / 2 = m - 1$\,.

(6)\,
Let us take any $\gp \in \assh{R}{R / \detid{2}{M}}$
and $n \geq m -1$\,.
Then, by (4) we have $\height{}{\gp} = m \geq 3$ and
$\gp \in \ass{R}{R / I^n}$\,.
Hence $\ass{R}{R / I^n}$ is {\it not} a subset of
$\{\,\gm\,\} \cup \Min{R}{R / I}$\,.
Therefore, by the observation stated in Introduction,
we get $\sat{I^n} \subsetneq I^{(n)}$
and the proof is complete.

\section{Computing $\sat{I^m}$}

In this section, we assume that
$R$ is an $(m+1)$-dimensional Cohen-Macaulay local ring
with the maximal ideal $\gm$ and
$x_1, x_2, \dots, x_{m+1}$ is an sop for $R$\,,
where $2 \leq m \in \zz$\,.
Let $M$ be the matrix stated in \ref{3.6}.
We put $I = \detid{m}{M}$\,.
Then, by (3) of \ref{3.6},
we have $\sat{I^n} = I^n$ for $1 \leq \forall n < m$\,.
The purpose of this section is to study $\sat{I^m}$\,.

For $1 \leq \forall j \leq m + 1$\,,
we set $a_j = (-1)^{j-1}\cdot\det M_j$\,,
where $M_j$ is the submatrix of $M$ determined by
removing the $j$-th column.
Then $I = (a_1, a_2, \dots , a_{m+1})R$\,.
Furthermore, for $1 \leq \forall k \leq m + 1$\,,
we denote by $\beta_k$ the minimum of the exponents of $x_k$
that appear in the entries of $M$\,.
Let us notice that $M$'s entries which are powers of $x_k$
appear as follows:
\[
\left(\begin{array}{cccc|c|cccc}
 & & & x_k^{\alpha_{1, k}} & \phantom{x^\alpha} & & & &  \\
 & & x_k^{\alpha_{2, k-1}} & & & & &  \\
 & \iddots & & & & & & &   \\
 x_k^{\alpha_{k, 1}} & & & & & & & & \phantom{x_{\beta_{\beta_\beta}}}  \\ \hline
 \phantom{x_k^{{\alpha_k}^\beta}}& & & & & & & & x_k^{\alpha_{k+1, m+1}} \\
 & & & & & & & x_k^{\alpha_{k+2, m}} &  \\
 & & & & & & \iddots & &   \\
 & & & & & x_k^{\alpha_{m, k+2}} & & &  
\end{array}\right)
\]
if $1 \leq k < m$\,, and 
\[
\left(\begin{array}{cccc|c}
 & & & x_m^{\alpha_{1, m}} & \phantom{x_m^\alpha} \\
 & & x_m^{\alpha_{2, m-1}} &   \\
 & \iddots & & &   \\
 x_m^{\alpha_{m, 1}} & & & &   \\
\end{array}\right)
\hspace{2ex}\mbox{or}\hspace{2ex}
\left(\begin{array}{c|cccc}
\phantom{x_m^\alpha} & & & & x_{m+1}^{\alpha_{1, m+1}} \\
 & & & x_{m+1}^{\alpha_{2, m}} \\
 & & \iddots & &  \\
 & x_{m+1}^{\alpha_{m, 2}} & & & 
\end{array}\right)
\]
if $k = m$ or $m + 1$\,, respectively.
So, we set
\[
\beta_k = \left\{\begin{array}{l}
\min\,\{ \alpha_{i, k-i+1} \}_{1 \leq i \leq k} \cup
     \{ \alpha_{i, k-i+m+2} \}_{k < i \leq m}
   \hspace{2ex}
   \mbox{if $1 \leq k < m$}\,,  \vspace{1em}\\
\min\,\{ \alpha_{i, k-i+1} \}_{1 \leq i \leq m}
   \hspace{2ex}
   \mbox{if $k = m$ or $m + 1$}\,.
\end{array}\right.
\]
Then, for $1 \leq \forall k \leq m + 1$\,,
we can choose $i_k$ with $1 \leq i_k \leq m$ so that
one of the following conditions is satisfied:
\[
\mbox{(i)\,
$1 \leq i_k \leq k$ and
$\beta_k = \alpha_{i_k, k-i_k+1}$
\hspace{2ex} or \hspace{2ex}
(ii)\,
$k < i_k \leq m$ and
$\beta_k = \alpha_{i_k, k-i_k+m+2}$\,.
}
\]
Now, for $1 \leq \forall i \leq m$ and
$1 \leq \forall k \leq {m + 1}$\,, we set
\[
x'_{ik} = \left\{\begin{array}{ll}
x_k^{\alpha_{i, k-i+1}-\beta_k} & \mbox{if $i \leq k$}\,, \vspace{1em}\\
x_k^{\alpha_{i, k-i+m+2}-\beta_k} & \mbox{if $i > k$}\,. 
\end{array}\right.
\]
Then $x'_{i_k, k} = 1$ for $1 \leq \forall k \leq m + 1$\,.
The next assertion can be verified easily.

\begin{lem}\label{4.1}
Suppose $1 \leq i \leq m$ and $1 \leq j \leq m + 1$\,.
\begin{itemize}
\item[{\rm (1)}]
If $i + j \leq m + 2$\,,
setting $k = i + j - 1$\,,
we have $1 \leq k \leq m + 1$\,,
$i \leq k$ and $x_{ij} = x_k^{\beta_k}\cdot x'_{ik}$\,.
\item[{\rm (2)}]
If $i + j > m + 2$\,,
setting $k = i + j - m - 2$\,,
we have $1 \leq k < m$\,,
$i > k$ and $x_{ij} = x_k^{\beta_k}\cdot x'_{ik}$\,.
\end{itemize}
\end{lem}

Let $Q$ be the ideal of $R$ generated by
$x_1^{\beta_1}, x_2^{\beta_2}, \dots , x_{m+1}^{\beta_{m+1}}$\,.
Then $M \in \setmat{m}{m+1}{Q}$ by \ref{4.1}.
The first main result of this section is the following:

\begin{thm}\label{4.2}
$\sat{I^m} = I^m :_R Q$ and
$\sat{I^m} / I^m \cong R / Q$\,.
\end{thm}

\noindent
{\it Proof}.\,
Let $S$ be the polynomial ring over $R$ with variables
$T_1, T_2, \dots , T_{m+1}$\,.
We regard $S$ as a graded ring by setting $\deg T_j = 1$
for $1 \leq \forall j \leq m + 1$\,. Let
\[
f_i = \sum_{j = 1}^{m+1} x_{ij}T_j \in S_1
\]
for $1 \leq \forall i \leq m$ and let
$K_\bullet$ be the Koszul complex of
$f_1, f_2, \dots , f_m$\,.
Then $K_\bullet$ is a graded complex.
Let $\partial_\bullet$ be the boundary map of $K_\bullet$
and let $e_1, e_2, \dots , e_m$ be an
$S$-free basis of $K_1$ consisting of homogeneous elements 
of degree $1$ such that
$\partial_1(e_i) = f_i$ for $1 \leq \forall i \leq m$\,.
As is stated in the proof of \ref{3.1},
\[
\begin{array}{ccccccccccccccc}
0 & \sra & [K_m]_m & \shomom{\partial_m} & 
  [K_{m-1}]_m & \sra & \cdots & \sra &
  [K_1]_m & \shomom{\partial_1} & [K_0]_m & 
  \shomom{\epsilon} & R & \sra & 0 \\
 & & & & & & & & & & \hspace{-0.5ex}\parallel & & & &  \\ 
 & & & & & & & & & & \,S_m & & & &
\end{array}
\]
is an acyclic complex,
where $\epsilon$ is the $R$-linear map such that
\[
\epsilon(T_1^{\alpha_1}T_2^{\alpha_2}\cdots T_{m+1}^{\alpha_{m+1}})
= a_1^{\alpha_1}a_2^{\alpha_2}\cdots a_{m+1}^{\alpha_{m+1}}
\]
for any $0 \leq \alpha_1, \alpha_2, \dots , \alpha_{m+1} \in \zz$
with $\alpha_1 + \alpha_2 + \cdots + \alpha_{m+1} = m$\,.
We obviously have ${\rm Im}\,\epsilon = I^m$\,.
We set $e = e_1 \wedge e_2 \wedge \cdots \wedge e_m$ and
$\check{e}_i = e_1 \wedge \cdots \wedge \widehat{e_i} \wedge
\cdots \wedge e_m$ for $1 \leq \forall i \leq m + 1$\,.
Let us take $\{\,e\,\}$ and 
$\{ T_j\check{e}_i \mid 1 \leq i \leq m\,,\,1 \leq j \leq m+1 \}$
as $R$-free basis of $[K_m]_m$ and $[K_{m-1}]_m$\,, respectively.
Because
\[
(\sharp) \hspace{6ex}
\partial_m(e) =
\sum_{i = 1}^m (-1)^{i-1} f_i\cdot\check{e}_i =
\sum_{i=1}^m \sum_{j = 1}^{m+1} (-1)^{i-1}x_{ij}\cdot T_j\check{e}_i\,,
\]
we have $\partial_m([K_m]_m) \subseteq Q\cdot [K_{m-1}]_m$\,.
Hence, by \cite[3.1]{i} we get
\[
(I^m :_R Q) / I^m \cong [K_m]_m / Q [K_m]_m \cong R / Q\,.
\]
Here, for $1 \leq \forall i \leq m$ and $1 \leq \forall k \leq m + 1$\,,
we set
\[
T_{ik} = \left\{
\begin{array}{ll}
T_{k-i+1} & \mbox{if $i \leq k$}\,, \vspace{1em} \\
T_{k-i+m+2} & \mbox{if $i > k$}\,.
\end{array}
\right.
\]
Then the following assertion holds:

\vspace{1em}
\noindent
{\bf Claim 1}\,
{\it
Suppose $1 \leq k, \ell \leq m + 1$ and $k \neq \ell$\,,
then $T_{ik} \neq T_{i\ell}$\,.
}

\vspace{1em}
\noindent
In order to prove the claim above,
we may assume $k < \ell$\,.
Then the following three cases can happen:
(i)\,
$i \leq k < \ell$
\,,\,
(ii)\,
$k \ < i \leq \ell$
\hspace{1ex}or\hspace{1ex}
(iii)\,
$k < \ell < i$\,.
Because $k - i + 1 < \ell - i + 1$ and
$k - i + m + 2 < \ell - i + m + 2$\,,
we get $T_{ik} \neq T_{i\ell}$ in the cases of (i) and (ii).
Furthermore, as $m + 1 > \ell - k$\,,
we get $k - i + m + 2 > \ell - i + 1$\,,
and so $T_{ik} \neq T_{i\ell}$ holds also in the case of (ii).
Thus we have seen Claim 1.

Now, for $1 \leq \forall k \leq m + 1$\,, we set
\[
v_{(k, e)} = \sum_{i = 1}^m (-1)^{i-1}x'_{ik}\cdot T_{ik}\check{e}_i\,.
\]
Then the following equality holds:

\vspace{1em}
\noindent
{\bf Claim 2}\,
$\partial_m(e) = \displaystyle{
\sum_{k = 1}^{m + 1} x_k^{\beta_k}\cdot v_{(k, e)}\,.}$

\vspace{1em}
\noindent
In fact, by ($\sharp$) and \ref{4.1} we have
\begin{eqnarray*}
\partial_m(e) & = &
  \sum_{i = 1}^m\,(\sum_{j = 1}^{m-i+2} (-1)^{i-1}x_{ij}\cdot T_j\check{e}_i +
               \sum_{j = m-i+3}^{m+1} (-1)^{i-1}x_{ij}\cdot T_j\check{e}_i\,)  \\
 & = &
  \sum_{i = 1}^m\,(\sum_{k = 1}^{m+1} (-1)^{i-1}x_k^{\beta_k}x'_{ik}\cdot
  T_{k-i+1}\check{e}_i +
  \sum_{k = 1}^{i-1} (-1)^{i-1}x_k^{\beta_k}x'_{ik}\cdot
  T_{k-i+m+2}\check{e}_i\,)   \\
 & = &
  \sum_{i = 1}^m\sum_{k = 1}^{m+1} (-1)^{i-1}x_k^{\beta_k}x'_{ik}\cdot
  T_{ik}\check{e}_i\,,
\end{eqnarray*}
and so the equality of Claim 2 follows.

Finally, we need the following:

\vspace{1em}
\noindent
{\bf Claim 3}\,
$\{ v_{(k, e)} \}_{1 \leq k \leq m}$ is a part of an $R$-free basis
of $[K_{m-1}]_m$\,.

\vspace{1em}
\noindent
If this is true, by \cite[3.4]{i}
(See \cite[3.4]{fin} for the case where $m = 2$)
we get $\dep{}{R / (I^m :_R Q)} > 0$\,,
which means $\sat{I^m} = I^m :_R Q$\,.
So, let us prove Claim 3.
By Claim 1, we see that
$T_{i_1, 1}\check{e}_{i_1}, T_{i_2, 2}\check{e}_{i_2}, \dots ,
T_{i_m, m}\check{e}_{i_m}$ are different to each other.
We set
\[
U = \{\,
T_j\check{e}_i \mid 1 \leq i \leq m\,,\,
1 \leq j \leq m + 1 \} \setminus
\{ T_{i_k, k}\check{e}_{i_k} \mid 1 \leq k \leq m \}
\]
and aim to prove that
$U \cup \{ v_{(k, e)} \}_{1 \leq k \leq m}$
is an $R$-free basis of $[K_{m-1}]_m$\,.
By \cite[3.3]{i}, it is enough to show that the submodule of
$[K_{m-1}]_m$ generated by
$U \cup \{ v_{(k, e)} \}_{1 \leq k \leq m}$
includes $T_{i_k, k}\check{e}_{i_k}$ for $1 \leq \forall k \leq m$\,.
This can be easily seen since
\[
v_{(k, e)} = (-1)^{k-1}\cdot T_{i_k, k}\check{e}_{i_k} +
\sum_{i \neq i_k} \,(-1)^{i-1}x'_{ik}\cdot T_{ik}\check{e}_i
\]
and $T_{ik}\check{e}_i \in U$ if $i \neq i_k$\,,
which follows from Claim 1.
Thus the assertion of Claim 3 follows,
and the proof of \ref{4.2} is complete. 

\vspace{1em}
If we assume a suitable condition on $\{ \alpha_{ij} \}$\,,
we can describe a generator of $\sat{I^m} / I^m$\,.
For $1 \leq \forall i \leq m$ and
$1 \leq \forall k \leq m + 1$\,, we set
\[
a_{ik} = \left\{\begin{array}{ll}
x'_{ik}a_{k-i+1} & \mbox{if $i \leq k$}\,, \vspace{0.5em} \\
x'_{ik}a_{k-i+m+2} & \mbox{if $i > k$}\,,
\end{array}\right.
\]
and $A = ( a_{ik} ) \in \setmat{m}{m+1}{I}$\,.
Then the next equality holds:

\begin{lem}\label{4.3}\hspace{1ex}
$\displaystyle{
A\,\left(
\begin{array}{c}
x_1^{\beta_1} \\
x_2^{\beta_2} \\
\vdots \\
x_{m+1}^{\beta_{m+1}}
\end{array}
\right)
=
\left(
\begin{array}{c}
0 \\
0 \\
\vdots \\
0
\end{array}\right)
}$\,.
\end{lem}

\noindent
{\it Proof}.\,
For $1 \leq \forall i \leq m$\,, we have
\[
\sum_{j = 1}^{m+1} x_{ij}a_j = 0\,.
\]
Let us divide the left side of this equality as follows:
\[
\sum_{j = 1}^{m-i+2} x_{ij}a_j + \sum_{j = m-i+3}^{m+1} x_{ij}a_j = 0\,.
\]
If $1 \leq j \leq m - i + 2$\,,
setting $k = i + j - 1$\,,
we have $i \leq k \leq m + 1$ and
\[
x_{ij}a_j = x_k^{\beta_k}x'_{ik}\cdot a_{k-i+1} = 
x_k^{\beta_k}\cdot a_{ik}\,.
\]
On the other hand, if $m - i + 3 \leq j \leq m + 1$\,,
setting $k = i + j - m - 2$\,,
we have $1 \leq k < i$ and
\[
x_{ij}a_j = x_k^{\beta_k}x'_{ik}\cdot a_{k-i+m+2} = 
x_k^{\beta_k}\cdot a_{ik}\,.
\]
Thus we get
\[
\sum_{k = 1}^{m+1} a_{ik}\cdot x_k^{\beta_k} = 0
\]
for $1 \leq \forall i \leq m$\,,
which means the required equality.

\vspace{1em}
For $1 \leq \forall k \leq m + 1$\,,
we denote by $A_k$ the submatrix of $A$
determined by removing the $k$-th column.
We set $b_k = \det A_k$\,.

\begin{ex}\label{4.4}
Suppose $\beta_k = \alpha_{k, 1}$ for $1 \leq \forall k \leq m$
{\rm (}For example, this holds if $\alpha_{k, 1} = 1$ for
$1 \leq \forall k \leq m${\rm )}.
Then, there exists $\delta \in R$ such that
$x_k^{\beta_k}\cdot\delta = b_k$ for
$1 \leq \forall k \leq m + 1$ and $\sat{I^m} = I^m + (\delta)$\,.
\end{ex}

\noindent
{\it Proof}.\,
The existence of $\delta$ such that
$x_k^{\beta_k}\cdot\delta = b_k$ for $1 \leq \forall k \leq m + 1$
follows from \ref{2.3} and \ref{4.3}.
Then $\delta \in I^m :_R Q \subseteq \sat{I^m}$\,.
We put $Q' = (x_1^{\beta_1}, x_2^{\beta_2}, \dots , x_m^{\beta_m})R$\,.
Then
\[
M_1 \equiv
\left(
\begin{array}{cccc}
   &    &    & x_{m+1}^{\,\alpha_{1, m+1}} \\
\multicolumn{2}{c}{\raisebox{0pt}[0pt][0pt]{\LARGE $0$}} 
            & x_{m+1}^{\,\alpha_{2, m}} &    \\
   & \iddots & \multicolumn{2}{c}{\raisebox{0pt}[0pt][0pt]{\LARGE $0$}} \\
x_{m+1}^{\,\alpha_{m, 2}} &   &   & 
\end{array}\right)
\,\mbox{mod}\,Q'\,,
\]
and so $a_1 \equiv \pm x_{m+1}^\alpha$ mod $Q'$\,,
where $\alpha := \alpha_{1, m+1} + \alpha_{2, m} + \cdots + \alpha_{m, 2}$\,.
Furthermore, if $2 \leq k \leq m + 1$\,,
we have $a_k \in Q'$ since the entries of the first column of $M_k$
are $x_1^{\beta_1}, x_2^{\beta_2}, \dots , x_m^{\beta_m}$\,.
Hence $Q' + I = Q' + (x_{m+1}^\alpha)$\,.
On the other hand, the assumption of \ref{4.4} implies that,
for $1 \leq \forall k \leq m$,
we can take $k$ itself as $i_k$\,,
and then $x'_{kk} = 1$\,.
Hence
\[
A = \left(\begin{array}{cccccc}
a_1 & x'_{12}a_2 & x'_{13}a_3 & \cdots & x'_{1m}a_m & x'_{1, m+1}a_{m+1} \\
x'_{21}a_{m+1} & a_1 & x'_{23}a_2 & \cdots & x'_{2m}a_{m-1} & x'_{2, m+1}a_m \\
x'_{31}a_m & x'_{32}a_{m+1} & a_1 & \cdots & x'_{3m}a_{m-2} & x'_{3, m+1}a_{m-1} \\
\hdotsfor{6} \\
x'_{m1}a_3 & x'_{m2}a_4 & x'_{m3}a_5 & \cdots & a_1 & x'_{m, m+1}a_2
\end{array}\right)\,,
\]
and so
\[
A_{m+1} \equiv
\left(\begin{array}{cccc}
\pm x_{m+1}^\alpha &   &   &     \\
  & \pm x_{m+1}^\alpha & 
\multicolumn{2}{c}{\raisebox{0pt}[0pt][0pt]{\LARGE $0$}} \\
\multicolumn{2}{c}{\raisebox{0pt}[0pt][0pt]{\LARGE $0$}} & \ddots &  \\
  &  &  & \pm x_{m+1}^\alpha
\end{array}\right)
\,\mbox{mod}\,Q'\,,
\]
which means
$b_{m+1} \equiv \pm x_{m+1}^{m\alpha}$ mod $Q'$\,.
Thus we get
\[
x_{m+1}^{\beta_{m+1}}\cdot\delta \equiv x_{m+1}^{m\alpha}
\,\mbox{mod}\,Q'\,.
\]
Here we notice $\beta_{m+1} \leq \alpha_{1, m+1} < \alpha$\,.
Because $x_1^{\beta_1}, x_2^{\beta_2}, \dots , x_{m+1}^{\beta_{m+1}}$
is an $R$-regular sequence, it follows that
\[
\delta \equiv \pm x_{m+1}^{m\alpha-\beta_{m+1}}
\,\mbox{mod}\,Q'\,,
\]
and so
\[
Q' + (\delta) = Q' + (x_{m+1}^{m\alpha-\beta_{m+1}}) \supseteq
Q' + (x_{m+1}^{m\alpha}) = Q' + I^m\,.
\]
Now we consider the $R$-linear map
\[
f : R \ra \frac{Q' + (x_{m+1}^{m\alpha-\beta_{m+1}})}{Q' + (x_{m+1}^{m\alpha})}
= \frac{Q' + (\delta)}{Q' + I^m}
\]
such that $f(1)$ is the class of $x_{m+1}^{m\alpha-\beta_{m+1}}$\,.
Then we have the following:

\vspace{1em}
\noindent
{\bf Claim}\hspace{1ex}
${\rm Ker}\,f = Q$\,.

\vspace{1em}
\noindent
If this is true,
then $R / Q \cong (Q' + (\delta)) / (Q' + I^m)$\,, and so
\[
\length{R}{R / Q} = \length{R}{\frac{Q' + (\delta)}{Q' + I^m}}\,.
\]
Because $(Q' + (\delta)) / (Q' + I^m)$ is a homomorphic image of
$(I^m + (\delta)) / I^m$ and
$I^m + (\delta) \subseteq \sat{I^m}$\,, we have
\[
\length{R}{\frac{Q' + (\delta)}{Q' + I^m}} \leq
\length{R}{\frac{I^m + (\delta)}{I^m}} \leq
\length{R}{\sat{I^m} / I^m} = \length{R}{R / Q}\,,
\]
where the last equality follows from \ref{4.2}\,.
Thus we see
\[
\length{R}{\frac{I^m + (\delta)}{I^m}} = \length{R}{\sat{I^m} / I^m}\,,
\]
and so $I^m + (\delta) = \sat{I^m}$ holds.

\vspace{1em}
\noindent
{\it Proof of Claim}.\,
Let us take any $r \in {\rm Ker}\,f$\,.
Then, there exists $s \in R$ such that
\[
r\cdot x_{m+1}^{m\alpha-\beta_{m+1}} \equiv s\cdot x_{m+1}^{m\alpha}
\,\mbox{mod}\,Q'\,.
\]
This congruence implies
\[
x_{m+1}^{m\alpha-\beta_{m+1}}(r - s\cdot x_{m+1}^{\beta_{m+1}}) \in Q'\,.
\]
Because $x_1^{\beta_1}, \dots , x_m^{\beta_m}, x_{m+1}^{m\alpha-\beta_{m+1}}$
is an $R$-regular sequence, we have
$r - s\cdot x_{m+1}^{\beta_{m+1}} \in Q'$\,,
which means $r \in Q$\,.
Hence ${\rm Ker}\,f \subseteq Q$\,.
As the converse inclusion is obvious,
we get the equality of the claim,
and the proof of \ref{4.4} is complete.

\end{document}